\documentclass[12pt,reqno]{amsart}
\makeatletter
\usepackage{amsthm}
\usepackage{amsmath}
\usepackage{version}
\usepackage{enumerate}
\usepackage{graphicx}
\usepackage[usenames,dvipsnames] {color}
\numberwithin{equation}{section}

\def\d{\mathbb{D}}
\def\c{\mathbb{C}}
\def\t{\mathbb{T}}

\def\cald{\mathcal{D}}
\def\calr{\mathcal{R}}

\def\calf{\mathcal{F}}
\newcommand\calm{\mathcal{M}}

\def\car#1{|#1|_{\rm car}}

\def\be{\begin{equation}}
\def\ee{\end{equation}}

\def\s0{s_0}
\def\p0{p_0}

\newcommand\nd{nondegenerate }

\DeclareMathOperator{\Car}{{\mathrm Car}}

 \newtheorem{theorem}{Theorem}[section]
 
 \newtheorem{lemma}[theorem]{Lemma}

\newtheorem{definition}[theorem]{Definition}

\newtheorem{defn}{Definition}[section]
\newtheorem{cor}[defn]{Corollary}

\newtheorem{rem}[defn]{Remark}

\newtheorem{fact*}{Fact}

\DeclareMathOperator\re{Re}

\DeclareMathOperator\rank{rank}
\DeclareMathOperator\hol{Hol}
\DeclareMathOperator\aut{Aut}

\newcommand\ov{\overline}
\newcommand\half{{\tfrac 12}}

\newcommand\dd{\mathrm d}
\newcommand\id{{\mathrm{id}}}
\newcommand\idd{\mathrm{id}_\mathbb{D}}

\newcommand\ess{\mathcal{S}}

\newcommand\lam{\lambda}

\newcommand{\T}{\mathbb{T}}

\newcommand{\X}{\mathcal{X}}
\newcommand{\D}{\mathbb{D}}
\newcommand{\C}{\mathbb{C}}

\newcommand{\ip}[2]{\left\langle #1, #2 \right\rangle}

\newcommand{\inv}{^{-1}}

\newcommand\fa{\mbox{ for all }}

\newcommand{\ph}{\varphi}
\renewcommand\phi{\varphi}
\newcommand\al{\alpha}

\newcommand\Ga{\Gamma}
\newcommand\de{\delta}

\newcommand\la{\lambda}

\newcommand\ups{\upsilon}

\newcommand\beq{\begin{equation}}
\newcommand\ds{\displaystyle}
\newcommand\eeq{\end{equation}}
\newcommand\df{\stackrel{\rm def}{=}}

\newcommand\black{\color{black}}

\newcommand\bbm{\begin{bmatrix}}
\newcommand\ebm{\end{bmatrix}}
\newcommand\bpm{\begin{pmatrix}}
\newcommand\epm{\end{pmatrix}}
\numberwithin{equation}{section}

\newtheorem{prop}[theorem]{Proposition}

\theoremstyle{definition}

\begin{document}
\title[Carath\'eodory extremal functions] {Nonuniqueness of Carath\'eodory extremal functions on the symmetrized bidisc}
%    author one information
\author{Jim Agler}
\address{Department of Mathematics, University of California at San Diego, CA \textup{92103}, USA}
%\curraddr{}
\email{jagler@ucsd.edu}
\thanks{Partially supported by National Science Foundation Grants
DMS 1361720 and 1665260 and by Newcastle University}

%    author two information
\author{Zinaida Lykova}
\address{School of Mathematics and Statistics, Newcastle University, Newcastle upon Tyne
 NE\textup{1} \textup{7}RU, U.K.}
%\curraddr{}
\email{Zinaida.Lykova@ncl.ac.uk}
\thanks{Partially supported by the Engineering and Physical Sciences grant EP/N03242X/1 and by the London Mathematical Society grant 42013}

%    author three information
\author{Nicholas Young}
\address{School of Mathematics and Statistics, Newcastle University, Newcastle upon Tyne NE1 7RU, U.K.
{\em and} School of Mathematics, Leeds University,  Leeds LS2 9JT, U.K.}
%\curraddr{}
\email{Nicholas.Young@ncl.ac.uk}
\date{22nd January, 2022}

\begin{abstract} 
We survey the Carath\'eodory extremal problem $\Car \de$ on the 
 symmetrized bidisc
\begin{eqnarray}
G &\stackrel{\rm{def}}{=}& \{(z+w,zw):|z|<1, \, |w|<1\} \notag\\
	&=& \{(s,p)\in \c^2: |s-\bar s p| < 1-|p|^2\}.  \notag
\end{eqnarray}
We also give some new results on this topic.
We are particularly interested in cases of this problem in which the solution of the problem is not unique. 
It is known that, for any $\de=(\lam,v)\in TG$ with $v\neq 0$, there is at least one $\omega\in\T$ such that $\Phi_\omega$  solves $\Car \de$, where $\Phi_\omega(s,p) = \frac{2\omega p-s}{2-\omega s}$.
Moreover, there is an essentially unique solution of $\Car \de$ if and only if $\de$ has exactly one Carath\'eodory extremal function of the form $\Phi_\omega$ for some $\omega\in\T$.
We give  a description of %all 
Carath\'eodory extremals for  $\de\in TG$ with more than one 
 Carath\'eodory extremal function  $\Phi_\omega$ for some values of $\omega \in\T$.  The proof exploits a model formula for the Schur class of $G$ which is an analog of the well-known network realization formula for Schur-class functions on the disc.
\end{abstract} 

\dedicatory{To Laszl\'o Lempert, with esteem and admiration}
\maketitle

\section*{Introduction}\label{intro}

In this paper we survey the theory  of the Carath\'eodory extremal problem $\Car \de$ in the symmetrized bidisc $G$. The problem is,  for a datum $\de= (\lam,v)$,  where  $\lam\in G$ and $v$ is a point in the holomorphic tangent space $T_\lam G\sim\c^2$ of $G$ at $\lam$, to understand the quantity
\beq\label{1}
\car{\de} \df   \sup_{F\in \hol(G,\d)} \frac{ |D_vF(\la)|}{1-|F(\la)|^2} 
\eeq
 and the functions $F$ for which the supremum in equation \eqref{1} is attained.
Here $\hol(G,\d)$ is the set of holomorphic maps from $G$ to the open unit disc $\d$ and $D_vF(\lam)$ is the directional derivative  of $F$ at $\lam$ in the direction $v$.
We call functions $F$ for which the supremum in equation \eqref{1} is attained {\em Carath\'eodory extremal functions } for the datum $\de=(\lam,v)$.

The Carath\'eodory extremal problem will forever be associated with the name of Lempert because of his astounding and beautiful theorem to the effect that the Carath\'eodory and Kobayashi metrics coincide on bounded convex domains \cite{lemp}.

 The domain we are interested in, the  {\em symmetrized bidisc}
\[
G \stackrel{\rm{def}}{=} \{(z+w,zw):|z|<1, \, |w|<1\},
\]
is a domain in $\c^2$ that has been the subject of many papers in the last 20 years, including \cite{jp04,cos04,ez05,sarkar,bhatta,aly2016,tryb} and many others.
Although $G$ is not convex (in fact it is not even isomorphic to any convex domain, see  \cite{cos04}), Lempert's conclusion of the equality of the Carath\'eodory and Kobayashi metrics nevertheless holds in $G$  \cite[Corollary 5.7]{ay2004}. 
In this paper we shall survey the Carathéodory extremal problem in $G$, and, in addition, shall
prove two new results on this problem in the case of non-unique extremal functions -- see Theorems \ref{simpleextremal} and \ref{pb}, which give sufficient and necessary conditions for the extremal functions in this case, but together fall short of an ``if and only if" statement.

Let us first establish some notation and describe the Carathéodory extremal problem in a general domain (that is, a connected open set) $\Omega$ in $\c^n$ and for a general datum $\de$ in $\Omega$.
For domains $\Omega_1,\ \Omega_2$, we denote by $\hol(\Omega_1,\Omega_2)$ the set of holomorphic maps from $\Omega_1 $ to $\Omega_2$.
By a {\em discrete datum} in $\Omega$ we mean an ordered pair $(\la,\mu)$ of points  of $\Omega$, while by an {\em infinitesimal datum} in $\Omega$ we mean an ordered pair  $\de=(\la,v)$ where $\la\in \Omega$ and $v$ is a point in the holomorphic tangent space $T_\la \Omega \sim \c^n$ of $\Omega$ at $\la$.  Thus an infinitesimal datum in $\Omega$ is the same thing as a point of the tangent bundle of $\Omega$.
 An infinitesimal datum in $\Omega$  will also be called a {\em tangent} (to $\Omega$).
 By a {\em datum} in $\Omega$ we mean either a discrete datum or an infinitesimal datum in $\Omega$.

   We say that $(\lam,v)$ is a {\em \nd} tangent if $v\neq 0$, and we say that a discrete datum $(\la,\mu)$ is \nd if $\la\neq\mu$.  We write $|\cdot|$ for the Poincar\'e metric on $T\d$:
\[
|(z, v)| \df   \frac{|v|}{1-|z|^2} \quad \mbox{ for } z\in\d, \mbox{ and } v\in\c,
\]
and $d$ for the Poincar\'e distance on $\d$:
\[
d(z,w) =  \tanh\inv \left| \frac{z-w}{1-\bar w z}\right|
\]
(we are following the conventions, though not the precise notation, of \cite[Chapter 1]{jp}).

The {\em Carath\'eodory}  or {\em Carath\'eodory-Reiffen  pseudometric} \cite{jp} on $\Omega$ is the Finsler pseudometric $\car{\cdot}$ on $T\Omega$ defined, for $\de=(\la,v)\in T\Omega$, by
\begin{align}\label{carprob}
\car{\de} & \df \sup_{F\in \hol(\Omega,\d)} |F_*(\de)| \notag \\
	&=   \sup_{F\in \hol(\Omega,\d)} \frac{ |D_vF(\la)|}{1-|F(\la)|^2}.
\end{align}
Here $F_*$ denotes the pushforward of $\de$ by the map $F$ to an element of $T\d$, given by
\[
\ip{g}{F_*(\de)}=\ip{g\circ F}{\de}
\]
for any analytic function $g$ in a neighbourhood of $F(\la)$.  Thus, if $\de=(\lam,v)$ then $F_*(\de)= (F(\lam), D_vF(\lam))$.

Analogously, the Carath\'eodory pseudodistance on $\Omega$ is the distance function given, for any discrete datum $\de=(\la,\mu)$, by
\begin{align}\label{carprob2}
\car{\de}  &=  \sup_{F\in \hol(\Omega,\d)} d(F(\la),F(\mu)) \notag \\
&=  \tanh\inv \black  \sup_{F\in \hol(\Omega,\d)}\left| \frac{F(\la)-F(\mu)}{1-\ov{F(\mu)}F(\la)}\right|
\end{align}

The {\em Carath\'eodory extremal problem}  $\Car \de$ on $\Omega$ is to study  the quantity $\car{\de}$  and to understand  the corresponding extremal functions.
We shall also say that  {\em $F$ solves the Carath\'eodory extremal problem}  $\Car \de$ on $\Omega$ if the supremum in either equation \eqref{carprob} or equation \eqref{carprob2} is attained at $F$.
  
  Observe that if $F$ is a Carath\'eodory extremal function for a datum $\de$ in $G$ then, by virtue of the invariance of the Poincar\'e distance and metric under automorphisms of $\d$, so is $m\circ F$ for any automorphism $m$ of $\d$.  This observation motivates us to say that   the solution of $\Car \de$ is {\em essentially unique } if, whenever $F_1$ and $F_2$ are Carath\'eodory extremal functions for $\de$, there exists $m\in\aut\d$ such that $F_2=m\circ F_1$.

  Consider any solution $F$ of $\Car\de$.  Among the functions $m\circ F$, with $m$ an automorphism of $\d$, there is exactly one that has the property
\[
m\circ F(\la)=0 \quad \mbox{ and } \quad D_v(m\circ F)(\la) > 0,
\]
or equivalently,
\beq\label{special}
(m\circ F)_*(\de) = (0,\car{\de}).
\eeq
Let us  say that a Carath\'eodory extremal function $F$ for a tangent  $\de=(\la,v)$ is {\em well aligned at} $\de$ if $F(\la)=0$ and  $D_v F(\la) > 0$ or, in other words, if the relation \eqref{special} holds.

It is rare, for a domain $\Omega$, that one can find either $\car{\cdot}$ or the corresponding extremal functions explicitly.  In some cases, however,  such as the ball in $\c^n$, the polydisc, and more generally, bounded symmetric homogeneous domains, there are more or less explicit formulae for $\car{\de}$.  An outlier is the case that $\Omega=G$, when the following result holds \cite[Theorem 1.1 and Corollary 4.3]{ay2004}.  We use the co-ordinates $\la=(s,p)$ for a point of $G$.  The unit circle will be denoted by $\T$.

\begin{theorem}\label{extthm10}
Let $\de$ be a \nd datum in  $G$.  There exists $\omega\in\t$ such that the function in $\hol(G,\d)$ given by
\beq\label{defPhi}
\Phi_\omega(s,p) \df \frac{2\omega p-s}{2-\omega s} \mbox{ for }(s,p)\in G
\eeq
solves $\Car \de$.
\end{theorem}

Theorem \ref{extthm10} implies that, for any datum $\de$ in $G$, the set 
\[
E_\de \df\{\omega\in\T: \Phi_\omega \mbox{ is a Carath\'eodory extremal for }\de \}
\]
 is non-empty.  In fact a ``trichotomy theorem" holds  (\cite[Theorem 1.6]{AY06}): for every datum $\de$ in $G$,  $E_\de$ either consists of a single point, or contains exactly $2$ points, or is equal to $\T$; moreover all three alternatives do occur. 
 
 In addition to the Carath\'eodory problem in $\Omega$, we shall also need the notion of a complex geodesic of a domain $\Omega$. 
 \begin{definition}\label{defcompgeo}
Let $\Omega$ be a domain and let $\cald\subset \Omega$.  We say that $\cald$ is a {\em complex geodesic} in $\Omega$ if there exists a function $k\in \hol(\d,\Omega)$ and a function $C\in \hol(\Omega,\d)$ such that $C\circ k=\id_\d$ and $\cald=k(\d)$.
\end{definition}

  Thus a complex geodesic is the image of an analytic disc $k:\d\to\Omega$ which is isometric with respect to the Poincar\'e distance $d$ on $\d$ and the Carath\'eodory distance $C_\Omega(\lam,\mu)=\car{(\lam,\mu)}$ on $\Omega$.

  In \cite[Theorem 3.6]{aly2016} we described five types of complex geodesic in $G$, namely purely unbalanced, exceptional, purely balanced, flat and royal, and showed that every complex geodesic in $G$ is of exactly one of these types. Furthermore,  it is known that each non-degenerate tangent $\de\in TG$ touches a unique complex geodesic $\cald$  \cite[Theorem A. 10]{aly2016}, and so we can define the type of a non-degenerate tangent $\de\in TG$ to be the type of the complex geodesic that touches $\de$. Therefore every \nd tangent in $TG$ is of exactly one of the five above-named types.
  
How can these five types of tangent be identified in terms of the corresponding Carath\'eodory extremal problems?  
  
For a non-degerate tangent $\de\in TG$ there is an essentially unique solution of $\Car \de $ if and only if $\de$ is either purely unbalanced or exceptional \cite[Theorem 2.1]{aly2019} and  \cite[Theorem 5.3]{kos}, which is to say that $\de$ has the property that there is a unique $\omega\in\t$ such that $\Phi_\omega$ solves $\Car \de$. L. Kosi\'{n}ski and W. Zwonek  commented in \cite{kos} that not much is known about the set of all Carath\'eodory extremal functions for a general \nd tangent $\de \in TG$. In this paper we consider the remaining types of tangent -- flat, royal and purely balanced tangents -- and will show that there are many Carath\'eodory extremal functions.
  We shall describe all Carath\'eodory extremal functions for royal tangents, and many extremal functions for flat tangents in Section \ref{royal}. We analyse the case of purely balanced tangents in Section \ref{purebal}.

The main tool we use in Section \ref{purebal} is a model formula for analytic functions from $G$ to the closed unit disc $\ov{\d}$ proved in \cite{AY2017} and stated below as Definition \ref{defGmodel} and Theorem \ref{modelGthm}.    We also use some basic properties of linear fractional transformations with operator coefficients.  For the convenience of readers we state the relevant facts in Section \ref{linfrac}.

  For the general theory of the Carath\'eodory extremal problem we refer to the monographs \cite{kob98, jp}.

\section{The model formula for the symmetrized bidisc}\label{Uniq}

Recall  the following model formula \cite[Definition 2.1]{AY2017}.
In this section we use the symbols $\la,\mu$ for points of $G$ and co-ordinates $\lam=(s,p)$.

\begin{definition}\label{defGmodel}
A $G$-\emph{model}
 for a function $\ph$ on $G$ is a triple $(\calm,T,u)$ where $\calm$ is a separable Hilbert space, $T$ is a contractive linear operator acting on $\calm$ and $u:G \to \calm$ is an analytic map such that, for all $\la,\mu\in G$,
\beq\label{modelform}
 1-\overline{\ph(\mu)}\ph(\la)= \ip{ (1-\mu_T^* \la_T) u(\la)}{u(\mu)}_\calm
\eeq
where, for $\la=(s,p)\in G$,
\[
\la_T \df (2pT-s)(2-sT)\inv.
\]

A $G$-model $(\calm,T,u)$ is {\em unitary} if $T$ is a unitary operator on $\calm$.
\end{definition}

For any domain $\Omega$ we define the {\em Schur class } $\ess(\Omega)$ to be the set of holomorphic maps from $\Omega$ to the closed unit disc $\d^-$.

\begin{theorem}\label{modelGthm}
Let $\ph$ be a function on $G$.  The following three statements are equivalent.
\begin{enumerate}
\item [\rm (1)]  $\ph\in\ess(G)$;
\item [\rm (2)] $\ph$ has a $G$-model;
\item [\rm (3)] $\ph$ has a unitary $G$-model $(\calm, T, u)$.
\end{enumerate}
\end{theorem}
For the proof see \cite[Theorem 2.2]{AY2017}.
From a $G$-model of a function $\ph\in\ess(G)$ one may easily proceed by means of a standard lurking isometry argument to a realization formula of the form
\[
\ph(\lam)=A+B\lam_T(1-D\lam_T)\inv C, \quad \mbox{ all } \lam\in G,
\]
for $\ph$, where $\begin{bmatrix} A&B\\C&D\end{bmatrix}$ is a contractive or unitary colligation on $\c\oplus\calm$.  However, for the present purpose it is convenient to work directly from the $G$-model rather than from the associated realization formula.

We also require a long-established fact about $G$ \cite{ay2004}, related to the coincidence of the Carath\'eodory and Kobayashi metrics on $TG$ \cite{jp}.
\begin{lemma}\label{k=c}
If $\de$ is a \nd tangent to $G$ and $F$ solves $\Car\de$ then there exists $k$ in $\hol(\d,G)$ such that $F\circ k=\idd$.  Moreover, if $\psi$ is any  solution of $\Car\de$ then $\psi\circ k$ is an automorphism of $\d$.
\end{lemma}

\
We shall need some minor  measure-theoretic technicalities, which are proved in \cite[Section 2]{aly2019}.

\begin{lemma}\label{baspos}
Let $Y$ be a set and let 
\[
A:\t\times Y\times Y \to \c
\]
be a map such that
\begin{enumerate}[\rm (1)]
\item $A(\cdot,z,w)$ is continuous on $\t$ for every $z,w \in Y$;
\item $A(\eta,\cdot,\cdot)$ is a positive kernel on $Y$ for every $\eta\in\t$.
\end{enumerate}
Let $\calm$ be a separable Hilbert space, let $T$ be a unitary operator on $\calm$ with spectral resolution
\[
T=\int_\t \eta \, \dd E(\eta)
\]
and let $v:Y\to\calm$ be a mapping.  
Let
\beq\label{defC}
C(z,w)= \int_\t A(\eta,z,w)\, \ip{\dd E(\eta)v(z)}{v(w)}
\eeq
for all $z,w \in Y$.  Then  $C$ is a positive kernel on $Y$.
\end{lemma}

\begin{lemma}\label{meas2}
For $i,j=1,2,$ let $a_{ij}:\t\to\c$ be a continuous map and let each $a_{ij}$ have only finitely many zeros in $\t$. 
Let $\nu_{ij}$ be a complex-valued Borel measure on $\t$ such that, for every Borel set $\tau$ in $\t$,
\[
\bbm \nu_{ij}(\tau) \ebm_{i,j=1}^2 \geq 0.
\]
Let $X$ be a Borel subset of $\t$ and suppose that
\[
 \bbm a_{ij}(\eta) \ebm_{i,j=1}^2\;  \mbox{ is positive and of rank } 2  \mbox{ for all } \eta \in X.
\]
Let
\[
C= \bbm c_{ij}\ebm_{i,j=1}^2
\]
where
\[
c_{ij}= \int_X a_{ij}(\eta)\,  \dd \nu_{ij}(\eta)   \qquad \mbox{ for }i,j=1,2.
\]
If  $\rank C \leq 1$
then either $c_{11}=0$ or $c_{22}=0$.

\end{lemma}

\section{Flat and royal tangents}\label{royal}

In this section we shall describe a large class of Carath\'eodory extremals for  flat and royal tangents.

Recall that a flat tangent has the form
\beq\label{aflatgeo}
\de=\left( (\beta+\bar\beta z,z), c(\bar\beta,1)\right)
\eeq
for some $z\in\d$ and $c\neq 0$, where $\beta\in\d$.  Such a tangent touches the `flat geodesic'
\[
F_\beta\df \{ (\beta+\bar\beta w, w):w\in\d\}.
\]
The description depends on a remarkable property of sets of the form $\calr\cup F_\beta, \ \beta\in\d$: they have the norm-preserving extension property in $G$ \cite[Theorem 10.1]{aly2016}.  That is, if $g$ is any bounded analytic function on the variety $\calr\cup F_\beta$, then there exists an analytic function $\tilde g$ on $G$ such that $g=\tilde g|\calr\cup F_\beta$ and the supremum norms of $g$ and $\tilde g$ coincide.  Indeed, the proof of \cite[Theorem 10.1]{aly2016} gives an explicit formula for one such $\tilde g$ in terms of a Herglotz-type integral.  Let us call the norm-preserving extension $\tilde g$ of $g$ constructed in \cite[Chapter 10]{aly2016} the {\em special extension} of $g$ to $G$.

\begin{lemma}\label{modflatgeo}
Let $\de=(\la,v)$ be the flat tangent $\la=(\beta +\bar\beta z,z), v= c(\bar\beta,1)$ for some $\beta,z\in\D$ and some $c\in\C$.
Then
\be\label{flatcar}
\car{\de} =\frac{|c|}{1-|z|^2}.
\ee
\end{lemma}
\begin{proof}
We shall show by a straightforward calculation that $|(\Phi_\omega)_*\de|$ is independent of $\omega$ in $\T$,
from which it follows by Theorem \ref{extthm10} that  $|(\Phi_\omega)_*\de|=\car{\de}$ for all $\omega\in\T$.
We have
\begin{align*}
D_v \Phi_\omega(\la) &=  \left[c\bar\beta \frac{\partial}{\partial s}\frac{2\omega p -s}{2-\omega s} + c\frac{\partial}{\partial p}\frac{2\omega p -s}{2-\omega s} \right]_{s=\beta+\bar\beta z, \ p=z} \\
	&= \left[ c\bar\beta \frac{(2-\omega s)(-1)-(2\omega p-s)(-\omega)}{(2-\omega s)^2} + c\frac{2\omega}{2-\omega s} \right]_{s=\beta+\bar\beta z, \ p=z}  \\
	&=\frac{2c(2\omega-\bar\beta-\omega^2\beta)}{(2-\omega\beta-\omega\bar\beta z)^2}.
\end{align*}
Thus
\begin{align*}
|(\Phi_\omega)_* \de| &=  \left[ \frac{|D_v\Phi_\omega(\la)|}{1-|\Phi_\omega(\la)|^2}\right]_{\la=(\beta+\bar\beta z,\ z)} \\
	&= \frac{|2c||2 -\bar\omega\bar\beta-\omega\beta|}{|2-\omega\beta-\omega\bar\beta z|^2\left(1-\left| \frac{2\omega z-\beta-\bar\beta z}{2-\omega\beta-\omega \bar\beta z}\right|^2 \right)} \\
	&= \frac{4|c||1-\re(\omega\beta)|}{4(1-\re(\omega\beta))(1-|z|^2)} \\
	&= \frac{|c|}{1-|z|^2}
\end{align*}
for all $\omega\in\T$.  Hence equation \eqref{flatcar} holds.
\end{proof}

\begin{theorem}\label{flat-car} \cite[Theorem 4.1]{aly2019}
Let $\de$ be the flat tangent
\beq\label{aflatgeo2}
\de=\left( (\beta+\bar\beta z,z), c(\bar\beta,1)\right)
\eeq
to $G$, where $\beta\in\d$ and $c\in\c\setminus\{0\}$.  Let $\zeta,\eta$ be the points in $\d$ such that
\[
(2\zeta,\zeta^2)=(\beta+\bar\beta\eta, \eta) \in  \calr\cap F_\beta
\]
and let $m$ be the unique automorphism of $\d$ such that  
\[
m_*((z,c))=(0,\car{\de}).
\]  

For every function $h\in\ess(\d)$ such that $h(\zeta)=m(\eta)$ the special extension $\tilde g$ to $G$ of the function
\beq\label{defg}
g: \calr\cup F_\beta\to \d, \quad (2w,w^2)\mapsto h(w), \quad (\beta+\bar\beta w,w) \mapsto m(w)
\eeq
for $w\in\d$ is a well-aligned Carath\'eodory extremal function for $\de$.
\end{theorem}

It is a simple calculation to show that $\calr\cap F_\beta$ consists of a single point.

A tangent $\de\in TG$ is said to be {\em royal} if it is tangent to the `royal variety', which is the variety $\calr \df\{(s,p) \in G: s^2=4p\}$ in $G$.  We have $\calr =\{(2z,z^2): z\in\D\}$, and so $\calr$ is a properly embedded submanifold of $G$, and one sees that a nondegenerate tangent $\de=(\la,v)\in TG$ is tangent to $\calr$ if and only if $\la=(2z,z^2)$ and $v=2c(1,z)$ for some $z\in\D$ and $c\in\C\setminus\{0\}$.   It is easy to calculate that, for a royal tangent $\de$, $(\Phi_\omega)_*(\de)$ is independent of $\omega\in\T$, and therefore $\Phi_\omega$ is a Carath\'eodory extremal for $\de$ for {\em every} $\omega\in\T$.
  In \cite[Theorem 3.1]{aly2019} we described {\em all} extremal functions for $\Car\de$ for royal tangents $\de$.   The statement is the following. 
  
\begin{theorem}\label{royalthm}
Let $\de\in TG$ be the royal tangent
\beq\label{roytang}
\de=\left((2z,z^2),2c(1,z)\right)
\eeq
for some $z\in\d$ and $c\in\c\setminus\{0\}$.  
\begin{enumerate}[\rm (1)]
\item $\car{\de}=\frac{|c|}{1-|z|^2}$;
\item   A function $F\in\hol(G,\d)$ solves $\Car\de$ if and only if there exists an automorphism $m$ of $\d$ and $\psi\in\ess(G)$ such that, for all points $(s,p)\in G$,
\beq
F(s,p)=m\left(  \half s +\tfrac 14 (s^2-4p)\frac{\psi(s,p)}{1-\half s\psi(s,p)}\right).
\eeq
\end{enumerate}
\end{theorem}

\section{Purely balanced tangents}\label{purebal}

There are tangents $\de\in TG$ such that $\Phi_\omega$ solves $\Car\de$ for exactly two values of $\omega$ in $\t$; we call them {\em purely balanced} tangents
 \cite[Section 3]{aly2016}.    In this section we study the general solution of $\Car\de$  for these tangents $\de=(\la,v)$.
 Such tangents can be described concretely as follows.  For any hyperbolic automorphism $m$ of $\d$ (that is, any automorphism that has two fixed points $\omega_1$ and $\omega_2$ in $\t$) let $h_m$ in $\hol(\d,G)$ be given by
\[
h_m(z)=(z+m(z),zm(z))
\]
for $z\in\d$.   A purely balanced tangent has the form
\beq\label{pbexpress}
\de=(h_m(z),c h_m'(z))
\eeq
for some hyperbolic automorphism $m$ of $\d$, some $z\in\d$ and some $c\in \c\setminus\{0\}$.  It is easy to see that, for $\omega\in\t$, the composition $\Phi_\omega\circ h_m$ is a rational inner function of degree at most $2$ and that the degree reduces to $1$ precisely when $\omega$ is either $\bar\omega_1$ or $\bar\omega_2$.  Thus, for these two values of $\omega$ (and these only), $\Phi_\omega\circ h_m$ is an automorphism of $\d$. It follows that $\Phi_\omega$ solves $\Car \de$ if and only if $\omega=\bar\omega_1$ or $\bar\omega_2$.

Observe that if $F, G$ are solutions of $\Car\de$ that are well aligned at $\de$ then any convex combination of
$F$ and $G$ also solves $\Car\de$.   Hence, if $\omega_1,\omega_2$ are the two values of $\omega$ in $\t$ for which $\Phi_\omega$ solves $\Car\de$, $\ups$ is any automorphism of $\d$ and $r\in[0,1]$ and   $m_j\in\aut\d$ is such that $m_j\circ \Phi_{\omega_j}$ is well aligned at $\de$,     then the function
\[
F= \ups \circ \left( r m_1\circ\Phi_{\omega_1}+ (1-r)m_2\circ\Phi_{\omega_2}\right)
\]
also solves $\Car\de$.  In fact it is easy to write down a much larger class of solutions of $\Car\de$, parametrised by the Schur class of $G$ and the interval $[0,1]$.
\begin{theorem}\label{simpleextremal}
Let $\de=(\la,v)$ be a non-degenerate purely balanced tangent to $G$ and let $\omega_1,\omega_2$ be the two points $\omega\in\T$ for which $\Phi_\omega$ solves $\Car\de$.  Let $m_j\in\aut\D$ be such that $m_j\circ\Phi_{\omega_j}$ is well aligned at $\de$ and let $\phi_j=m_j\circ\Phi_{\omega_j}$ for $j=1,2$.
For any $\psi\in\ess(G)$ and $r\in [0,1]$ the function
\be\label{linfracexfn}
F\df r\ph_1+(1-r)\ph_2 + r(1-r) \frac {(\ph_1-\ph_2)^2\psi}{1-[(1-r)\ph_1+r\ph_2]\psi}
\ee
is a Carath\'eodory extremal for $\de$ and is well aligned at $\de$.
\end{theorem}
\begin{proof}
For $r\in [0,1]$ the matrix
\[
U_r=\bpm \sqrt{r} & \sqrt{1-r}\\  \sqrt{1-r} & -\sqrt{r} \epm
\]
is unitary.  Hence the matrix function
\[
T\df \bpm A & B\\C & D \epm \df U_r^*\bpm \ph_1 & 0\\0 & \ph_2 \epm U_r = \bpm r\ph_1+(1-r)\ph_2 & \sqrt{r(1-r)}(\ph_1-\ph_2) \\ \sqrt{r(1-r)}(\ph_1-\ph_2) & (1-r)\ph_1+r\ph_2 \epm
\]
satisfies $\|T(s,p)\|< 1$ for all  $(s,p)\in G$.  It follows from properties of linear fractional transformations (Corollary \ref{calFleq1} and Remark \ref{invertpartdefinedrealform}) that, for any $\psi\in\ess(G)$,  the function
\begin{align}\label{linearfractional}
F\df \calf_T\circ \psi &\df A+B\psi(1-D\psi)\inv C  \notag\\
	&=r\ph_1+(1-r)\ph_2 + r(1-r) \frac {(\ph_1-\ph_2)^2\psi}{1-[(1-r)\ph_1+r\ph_2]\psi}
\end{align}
is analytic on $G$ and maps $G$ to $\D$.  Since $\ph_1,\ph_2$ are well aligned at $\de$ we have $\ph_1(\la)=\ph_2(\la)=0$, so that $F(\la)=0$, and one finds that
\[
D_vF(\la) = r D_v\ph_1(\la) + (1-r)D_v\ph_2(\la) = \car{\de}.
\]
Thus $F$  is a holomorphic map from $G$ to $\D$ such that $F(\la)=0$ and $D_vF(\la)= \car{\de}$, that is, $F$ solves $\Car \de$ and is well aligned at $\de$.
\end{proof}

\begin{lemma}\label{estu}
Let $(\calm,T,u)$ be a unitary $G$-model for a function in $\ess(G)$.  For any $(s,p)\in G$,
\beq\label{1.1}
\|u(s,p)\| \leq \frac{4-|s|^2}{\left( (4-|s|^2)^2-(2|s-\overline{s}p|+|s^2-4p|)^2   \right)^{\half}}.
\eeq
\end{lemma}
\begin{proof}
Let $(\calm,T,u)$ be a unitary $G$-model for a function $F\in \ess(G)$ and let $T$ have spectral resolution
\[
T=\int_\t \eta \,  \dd E(\eta).
\]
Then the model formula for $F$ is
\[
1-\overline{F(\mu)}F(\la)=\int_\t 1-\overline{\Phi_\eta(\mu)}\Phi_\eta(\la) \, \ip{\dd E(\eta)u(\la)}{u(\mu)}
\]
for all $\la,\mu\in G$.  Therefore, for any $\la\in G$,
\begin{align*}
1 &\geq 1- |F(\la)|^2 \\
	&= \int_\t (1-|\Phi_\eta(\la)|^2) \, \ip{\dd E(\eta)u(\la)}{u(\la)}\\
	&\geq \inf_\eta (1-|\Phi_\eta(\la)|^2)  \ip{E(\t)u(\la)}{u(\la)} \\
	&= (1-\sup_\eta|\Phi_\eta(\la)|^2) \|u(\la)\|^2.
\end{align*}
For any $(s,p)\in G$, the linear fractional transformation $\eta \mapsto \Phi_\eta(s,p)$ maps $\t$ to the circle with centre and radius
\[
2\frac{s-\overline{s}p}{4-|s|^2} \mbox{ and } \frac{|s^2-4p|}{4-|s|^2}
\]
respectively.  Hence, if $\la=(s,p)\in G$, then
\[
\sup_\eta|\Phi_\eta(\la)|  =  \frac{2|s-\overline{s}p|+|s^2-4p|}{4-|s|^2}.
\]
Inequality \eqref{1.1} follows.
\end{proof}

It is conceivable that the converse of Theorem \ref{simpleextremal} is also true: every Carath\'eodory extremal function for a purely balanced tangent $\de$
can be written in the form \eqref{linearfractional} for some $\psi$ in the Schur class $\ess(G)$ and some $r\in (0,1)$.   We do not know whether this is so, but we do have another description of the general Carath\'eodory extremal function for a purely balanced tangent $\de\in TG$.

\begin{theorem}\label{pb}
Let $\de=(\la,v)$ be a purely balanced tangent to $G$ and let $\omega_1,\omega_2$ be the two points $\omega\in\T$ for which $\Phi_\omega$ solves $\Car\de$.  Let $m_j\in\aut\D$ be such that $m_j\circ\Phi_{\omega_j}$ is well aligned at $\de$ and let $\phi_j=m_j\circ\Phi_{\omega_j}$ for $j=1,2$.    Suppose that the automorphism $m_j$ is given by the formula
\be\label{formmj}
m_j(z) = c_j \frac{z-\al_j}{1-\ov{\al_j}z} \quad  \fa z\in\D,
\ee
where $|c_j|=1$ and $\al_j\in\D$.
For any well-aligned Carath\'eodory extremal function $F$ for $\de$ there exist holomorphic functions $u_1,u_2$ on $G$ satisfying, for all $\mu=(s,p)\in G$,
\be\label{ujsum1}
  \sum_{j=1}^2 \frac{(1-|\al_j|^2)^\half u_j(\mu)}{1+\bar\al_j \bar c_j\ph_j(\mu)}= 1
\ee
and
\be\label{estubis}
|u_1(\mu)|^2+|u_2(\mu)|^2 \leq    \frac{(4-|s|^2)^2}{(4-|s|^2)^2-(2|s-\bar s p|+|s^2-4p|)^2}
\ee
such that, for all $\mu\in G$,
\be\label{gotphi}
F(\mu)=  \sum_{j=1}^2 \frac{(1-|\al_j|^2)^\half \ph_j(\mu)u_j(\mu)}{1+\bar\al_j \bar c_j\ph_j(\mu)}.
\ee

\end{theorem}

\begin{proof}
Let $F$ be a Carath\'eodory extremal function for the purely balanced tangent $\de \in TG$ and suppose that $F$ is well aligned at $\de$.
By Lemma \ref{k=c} and the tautness of $G$ (see for example \cite{jp}), there exists a complex geodesic $k\in \hol(\d,G)$ such that $\de$ touches $k(\d)$, and so $\de=k_*(\de_0)$ for some tangent $\de_0\in T\d$.  Necessarily $\car{\de_0}=\car{\de}$.  On replacing $k$ by $k\circ \ups$, for a suitably chosen automorphism $\ups$ of $\d$, we can arrange that $\de_0=(0,\car{\de})$.  Then $F \circ k$ is an automorphism of $\d$ such that
\[
(F\circ k)_*(0,\car{\de})=F_*(\de)=(0,\car{\de}),
\]
and hence, by the Schwarz Lemma, 
\beq\label{phcirck}
F \circ k= \idd.
\eeq
By the choice of $m_j$,
\[
(m_j\circ\Phi_{\omega_j}\circ k)_*(0,\car{\de})= (m_j\circ\Phi_{\omega_j})_*(\de)=(0,\car{\de}),
\]
and so, again by the Schwarz Lemma,
\be\label{phojmj}
\Phi_{\omega_j}\circ k= m_j\inv, \    \quad \mbox{ for }j=1,2.
\ee
By choice of $m_j$,
\begin{align}\label{defphj}
\ph_j \df m_j\circ \Phi_{\omega_j}
\end{align}
 is well aligned at $\de$ for $ j=1,2$.

Choose a unitary $G$-model $(\calm,T,u)$ of $F$ and let
$E$ be the spectral measure of the unitary operator $T$, so that $T=\int_{\t} \eta \  \dd E(\eta)$.  Thus $\calm$ is a separable Hilbert space, $u:G \to \calm$ is a holomorphic map, and,  by the spectral theorem, for $\la=(s,p)\in G$,
\[
\la_T= (2pT-s)(2-sT)\inv = \int_\T \Phi_\eta(\la) dE(\eta),
\]
and so equation \eqref{modelform} yields the representation, for all $\la,\mu$ in $G$,
\beq\label{modelbis}
1-\overline{F(\mu)}F(\la) = \int_\t \left( 1- \overline{\Phi_\eta(\mu)}\Phi_\eta(\la)\right)\ip{\dd E(\eta)u(\la)}{u(\mu)}.
\eeq
Consider $z,w\in\d$, put $\la=k(z), \ \mu=k(w)$ in equation \eqref{modelbis} and divide by $1-\bar w z$ to obtain, for all $z,w\in\d$,
\begin{align}
1&= \frac{1-\overline{F \circ k(w)} F \circ k(z)}{1-\bar w z} \notag \\
	&=\int_{\{\omega_1,\omega_2\}} + \int_{\t\setminus\{\omega_1,\omega_2\}}  \frac{1-\overline{\Phi_\eta \circ k(w)}\Phi_\eta\circ k(z)}{1-\bar w z} \ip{\dd E(\eta)u\circ k(z)}{u\circ k(w)} \notag \\
	&= I_1+I_2 \label{page65}
\end{align}
where
\begin{align}
I_1(z,w)&= \int_{\{\omega_1,\omega_2\}}  \frac{1-\overline{\Phi_\eta \circ k(w)}\Phi_\eta\circ k(z)}{1-\bar w z} \ip{\dd E(\eta)u\circ k(z)}{u\circ k(w)}\notag\\	&=  \sum_{j=1}^2\frac{1-\ov{\Phi_{\omega_j}\circ k(w)} \Phi_{\omega_j}\circ k(z)}{1-\bar w z}\ip{E(\{\omega_j\})u\circ k(z)}{u\circ k(w)} \notag\\
	&= \sum_{j=1}^2 \frac{1-\overline{m_j\inv(w)}m_j\inv(z)}{1-\bar w z}\ip{E(\{\omega_j\}) u\circ k(z)}{u\circ k(w)}, \mbox{ by equation \eqref{phojmj}, }  \label{I1}\\
I_2(z,w) &= \int_{\t\setminus\{\omega_1,\omega_2\}}  \frac{1-\overline{\Phi_\eta \circ k(w)}\Phi_\eta\circ k(z)}{1-\bar w z} \ip{\dd E(\eta)u\circ k(z)}{u\circ k(w)}. \label{I2}
\end{align}
The left-hand side of equation \eqref{page65} is a positive kernel of rank $1$ on $\d$ and $I_1$ is also a positive kernel on $\d$.  The integrand in $I_2$ is a positive kernel on $\d$ for each $\eta\in\t$, by Pick's theorem, since $\Phi_\eta\circ k$ is in the Schur class.   Hence, by Lemma \ref{baspos}, $I_2$ is also a positive kernel on $\d$.  Since $I_1 +I_2$ has rank $1$, it follows that $I_2$ has rank at most $1$ as a kernel on $\d$.

We now modify the argument in the proof of  \cite[Theorem 2.1]{aly2019} to show that $E(\t\setminus \{\omega_1,\omega_2\}) u \circ k(w)=0$ for all $w\in\d$.   Since $\Phi_\eta$ does not solve $\Car\de$ for $\eta\in\t\setminus \{\omega_1,\omega_2\}$, for such $\eta$, $\Phi_\eta\circ k$ is a Blaschke product of degree $2$.  Hence, for any choice of distinct points $z_1,z_2$ in $\d$ and any $\eta\in\t\setminus \{\omega_1,\omega_2\}$, the $2\times 2$ matrix 
\[
\bbm a_{ij}(\eta)\ebm_{i,j=1}^2 \df \bbm \ds \frac{1-\overline{\Phi_\eta \circ k(z_i)}\Phi_\eta\circ k(z_j)}{1-\bar z_i z_j} \ebm_{i,j=1}^2
\]
 is positive and has rank $2$.  Each $a_{ij}$ is a ratio of nonzero trigonometric polynomials in $\eta$, hence has only finitely many zeros in $\t$.

Define Borel measures $\nu_{ij}, \ i,j=1,2$, on $\t\setminus \{\omega_1,\omega_2\}$ by
\[
\nu_{ij}=\ip{E(\cdot)u\circ k(z_i)}{u\circ k(z_j)}.
\]
Then $\bbm \nu_{ij}(\tau)\ebm\geq 0$ for every Borel subset $\tau$ of $ \t\setminus \{\omega_1,\omega_2\}$.
By equation \eqref{I2},
\[
I_2(z_i,z_j) = \int_{\t\setminus \{\omega_1,\omega_2\}} a_{ij} \, \dd \nu_{ij}.
\]
Moreover, by equation \eqref{page65},
\[
[I_2(z_i,z_j)] \leq [I_1(z_i,z_j)] + [I_2(z_i,z_j)] = \bbm 1&1\\1&1\ebm.
\]
It follows that 
\[
\bbm\int_{\t\setminus \{\omega_1,\omega_2\}}a_{ij}\, \dd\nu_{ij}\ebm = \bbm I_2(z_i,z_j)\ebm = \kappa \bbm 1&1\\1&1\ebm
\]
for some $\kappa\in [0,1]$. We may apply  Lemma \ref{meas2} with $X = \t\setminus \{\omega_1,\omega_2\}$ to deduce that $\kappa=0$. Hence $  [I_2(z_i,z_j)] =0$. In particular
\[ 0 = I_2(z_1,z_1) = \int_{\t\setminus \{\omega_1,\omega_2\}} a_{11} \, \dd \nu_{11}.
\]
Thus
\[  E(\t\setminus \{\omega_1,\omega_2\}) u \circ k(z_1)=0. \]
Since $z_1,z_2$ are arbitrary distinct points in $\d$,  $I_2\equiv 0$ and so, by equation \eqref{page65}, for all $z,w\in\d$,
\begin{align}
I_1(z,w)&= 1, \label{I1=1} \\
\ip{E(\t\setminus\{\omega_1,\omega_2\}) u\circ k(z)}{u\circ k(w)}&=I_2(z,w)= 0, \label{4.12}
\end{align}
wherefore
\[
E(\t\setminus\{\omega_1,\omega_2\}) u\circ k(z) = 0 \qquad\mbox{ for all }z\in \d.
\]

For $\al\in\d$ let $K_\al$ denote the normalised Szeg\H{o} kernel at $\al$:
\[
K_\al(z)= \frac{(1-|\al|^2)^{\half}}{1-\bar\al z} \quad \mbox{ for } z\in\d.
\]
For $j=1,2$, since the automorphism $m_j$ is given by equation \eqref{formmj},
we have, for $z\in\D$,
\[
m_j\inv(z) =  \bar c_j\frac{z+c_j\al_j}{1+\bar c_j \bar\al_j z},
\]
and so
\[
1-\ov{m_j\inv(w)}m_j\inv(z) = \frac{(1-|\al_j|^2)(1-\bar w z)}{(1+c_j\al_j\bar w)(1+\bar c_j\bar\al_j z)},
\]
and therefore
\[
\frac{1-\ov{m_j\inv(w)}m_j\inv(z)}{1- \bar w z} = \ov{K_{-c_j\al_j}(w)} K_{-c_j\al_j}(z).
\]

  We can thus write equations \eqref{I1} and \eqref{I1=1} in the form
\begin{align*}
1&= \sum_{j=1}^2   \overline{K_{-c_j\al_j}(w)}K_{-c_j\al_j}(z)  \ip{E(\{\omega_j\}) u\circ k(z)}{u\circ k(w)} \\
	&= \ip{ \bpm K_{-c_1\al_1}(z) E(\{\omega_1\}) u\circ k(z) \\ K_{-c_2\al_2}(z) E(\{\omega_2\}) u\circ k(z)\epm}{ \bpm K_{-c_1\al_1}(w) E(\{\omega_1\}) u\circ k(w) \\ K_{-c_2\al_2}(w) E(\{\omega_2\}) u\circ k(w) \epm}_{\calm \oplus\calm}.
\end{align*}
Now if $f$ is a Hilbert-space-valued function on $\D$ such that $\ip{f(z)}{f(w)} =1$ for all $z,w$, then it is immediate from consideration of $\|f(z)-f(w)\|^2$ that $f$ is constant, with value of unit norm.
Hence there exists a unit vector $x=(x_1,x_2) \in \calm\oplus\calm$ such that 
\beq\label{proptauj}
K_{-c_j\al_j}(z) E(\{\omega_j\}) u\circ k(z) = x_j \qquad \mbox{ for } j=1,2 \, \mbox{ and } z\in\d.
\eeq
Moreover, since $E(\{\omega_1\})$ and $E(\{\omega_2\})$ have orthogonal ranges, $x_1\perp x_2$ and $\|x_1\|^2 + \|x_2\|^2 =1$.

In the model formula \eqref{modelbis} choose $t=k(w)$ for some $w\in\d$ to obtain, for any $\la\in G$,
\begin{align}
1-\bar wF(\la) &= 1- \overline{F \circ k(w)}F(\la)\notag  \\
	&= \int_{\{\omega_1,\omega_2\}} + \int_{\t\setminus\{\omega_1,\omega_2\}} \left(1-\overline {\Phi_\eta\circ k(w)}\Phi_\eta(\la)\right) \ip{\dd E(\eta)u(\la)}{u\circ k(w)}\notag\\
	&= \sum_{j=1}^2 \left(1-\overline{m_j\inv(w)}\Phi_{\omega_j}(\la)\right) \ip{E(\{\omega_j\})u(\la)}{u\circ k(w)} + I_3 \label{powereq}
\end{align}
where
\beq\label{I2bis}
I_3(\la,w) =  \int_{\t\setminus\{\omega_1,\omega_2\}} \left(1-\overline {\Phi_\eta\circ k(w)}\Phi_\eta(\la)\right) \ip{\dd E(\eta)u(\la)}{u\circ k(w)}.
\eeq
Since $E(\t\setminus\{\omega_1,\omega_2\}) u\circ k(z) =0$ for all $z$, it follows that $I_3=0$ and so, for every $w\in\d$ and $\la\in G$,
 by equations \eqref{proptauj} and  \eqref{powereq},
\begin{align}
1-\bar w F(\la) &= \sum_{j=1}^2 (1-\overline{m_j\inv (w)}\Phi_{\omega_j}(\la)) \ip{u(\la)}{E(\{\omega_j\})u\circ k(w)}\notag\\
	&=  \sum_{j=1}^2 (1-\overline{m_j\inv(w)}\Phi_{\omega_j}(\la)) \ip{u(\la)}{x_j/K_{-c_j \al_j}(w)} \notag \\
	&=    \sum_{j=1}^2 \frac{1+c_j\al_j\bar w}{(1-|\al_j|^2)^{\half}} (1-\overline{m_j\inv(w)}\Phi_{\omega_j}(\la)) \ip{u(\la)}{x_j} \notag \\
	&=  \sum_{j=1}^2   \left(1+c_j\al_j \bar w - ( c_j \bar w+\bar\al_j)\Phi_{\omega_j}(\la)\right) \frac{\ip{u(\la)}{x_j}}{(1-|\al_j|^2)^{\half}}.
\label{gotit}
\end{align}

Consider first the case that $x_2=0$.  Then $\|x_1\|=1$ and equation \eqref{gotit} becomes
\[
1-\bar w F(\la)\ =\  \left(1+c_1\al_1 \bar w - (c_1 \bar w+\bar\al_1)\Phi_{\omega_1}(\la)\right) \frac{ \ip{u(\la)}{x_1}}{(1-|\al_1|^2)^{\half}}
\]
for all $w\in\D$ and $\la\in G$.  Put $w=0$ to obtain
\[
\frac{\ip{u(\la)}{x_1}}{(1-|\al_1|^2)^{\half}}\ =\ \frac{1}{1 - \ov{\al_1}\Phi_{\omega_1}(\la)}
\]
and then equate coefficients of $\bar w$ to deduce that
\begin{align*}
F(\la)\ &=\ c_1\left(\Phi_{\omega_1}(\la) - \al_1\right) \frac{\ip{u(\la)}{x_1}}{(1-|\al_1|^2)^{\half}} \\
	&=\  c_1\frac{\Phi_{\omega_1}(\la) - \al_1}{1-\ov{\al_1}\Phi_{\omega_1}(\la)}\\
         &=m_1\circ\Phi_{\omega_1}(\la)\\
	&=\ph_1(\la).
\end{align*}
Thus, in the case $x_2=0$,  $F$ is given by the formula \eqref{gotphi} with the choice 
$$u_1=\frac{1+\ov{\al_1}\ov{c_1}\ph_1}{(1-|\al_1|^2)^{\half}}$$ and $u_2=0.$
Similarly, in the case that $x_1=0$, the extremal function $F$ is $m_2\circ\Phi_{\omega_2}=\ph_2$, and $F$ corresponds to the choice $u_1=0$  in equation \eqref{gotphi}.

It remains to consider the case that $x_1$ and $x_2$ are both non-zero vectors.
Let 
\[
u_j(\la) =  \ip{u(\la)}{x_j} \quad \mbox{ for } j=1,2, \mbox{ and for all } \la\in G
\]
  Then $u_j\in\hol(G,\C)$ and since $u_1,u_2$ are the components of $u$ with respect to the orthonormal vectors $x_1,x_2$ in $\calm$,
 Lemma \ref{estu} implies that, for $\la=(s,p)\in G$,
\[
|u_1(\la)|^2+ |u_2(\la)|^2 \leq \|u(s,p)\|^2 \leq  \frac{(4-|s|^2)^2}{(4-|s|^2)^2-(2|s-\bar s p|+|s^2-4p|)^2},
\]
in agreement with the estimate \eqref{estubis}.
Equation \eqref{gotit} can be written
\[
1-\bar wF(\la)=\sum_{j=1}^2   \left(1-\bar \al_j\Phi_{\omega_j}(\la)- \bar w( c_j\Phi_{\omega_j}(\la) - c_j\al_j)\right)\frac{u_j(\la)}{(1-|\al_j|^2)^\half}.
\]
On equating the coefficients of $\bar w$  and $\bar w^0$ in  equation \eqref{gotit} we obtain the relations
\begin{align}\label{1phi}
1&= \sum_{j=1}^2 \left(1- \bar \al_j \Phi_{\omega_j}(\la)\right) \frac{u_j(\la)}{(1-|\al_j|^2)^\half},  \notag \\
F(\la) &=  \sum_{j=1}^2\left( c_j\Phi_{\omega_j}(\la) -c_j \al_j\right) \frac{u_j(\la)}{(1-|\al_j|^2)^\half}.
\end{align}
Since 
\[
\Phi_{\omega_j} = m_j\inv\circ \ph_j = \frac{\bar c_j\ph_j + \al_j}{1+\bar\al_j \bar c_j\ph_j},
\]
we may re-write equations \eqref{1phi} in the form
\begin{align}\label{1phia}
1&= \sum_{j=1}^2  \frac{(1-|\al_j|^2)^\half u_j(\la)}{1+\ov{\al_jc_j}\ph_j(\la)},  \notag \\
F(\la) &=  \sum_{j=1}^2 \frac{(1-|\al_j|^2)^\half \ph_j(\la)u_j(\la)}{1+\ov{\al_jc_j}\ph_j(\la)}
\end{align}
for all $\la \in G$, in agreement with Theorem \ref{pb}.
\end{proof}

\section{Contractive linear fractional transformations}\label{linfrac} 
This section is an appendix in which we recall an identity about linear fractional transformations of operators.  We presume that
this identity has long been known.  A special case of it was stated by C. L. Siegel in \cite{siegel}, and Siegel may well have known 
the identity in its general form.  We stated the identity without proof in previous publications, see \cite{realizsymmbidisc}. We were unable to find a proof of the identity in the literature, therefore, we refer the reader to a proof of it in
\cite[Section B.2]{brown}.  The proof is by straightforward algebraic calculation, but is not entirely trivial to work out from scratch.

Let  $H,G,U$ and $V$ be Hilbert spaces. Let $P$ be an operator such that 
\[P=\begin{bmatrix} P_{11} & P_{12} \\ P_{21} & P_{22} \end{bmatrix}:H\oplus U\to G\oplus V\]
and let $X:V\to U$ be an operator for which $I-P_{22}X$ is invertible in the space $\mathcal{B}(V)$ of bounded linear operators on $V$ . Then
we denote by $\mathcal{F}_P(X)$\index{$\mathcal{F}_P(X)$} the linear fractional transformation 
\[\mathcal{F}_P(X):=P_{11}+P_{12}X(I-P_{22}X)^{-1}P_{21}:H\to G.\] 

\begin{prop}\label{calF}\textup{\cite[Lemma 1.7]{realizsymmbidisc}} 
Let $H,G,U$ and $V$ be Hilbert spaces. Let \[P=\begin{bmatrix} P_{11} & P_{12} \\ P_{21} & P_{22} \end{bmatrix}\text{ and }
Q=\begin{bmatrix} Q_{11} & Q_{12} \\ Q_{21} & Q_{22} \end{bmatrix}\] be operators from $H\oplus U$ to $G\oplus V$, and
let $X$ and $Y$ be operators from $V$ to $U$ for which $I-P_{22}X$ and $I-Q_{22}Y$ are invertible in $\mathcal{B}(V)$.
Then
\begin{align*}I_H-\mathcal{F}_Q(Y)^*&\mathcal{F}_P(X)=\\
&=Q^*_{21}(I_V-Y^*Q^*_{22})^{-1}(I_V-Y^*X)(I_V-P_{22}X)^{-1}P_{21}\\
&\qquad+\begin{bmatrix}I_H& Q^*_{21}(I_V-Y^*Q^*_{22})^{-1}Y^* \end{bmatrix}(I_{H\oplus U}-Q^*P)
\begin{bmatrix} I_H\\X(I_V-P_{22}X)^{-1}P_{21} \end{bmatrix}.\end{align*}
\end{prop}

\begin{cor}\label{calFleq1}
Let $H,G,U$ and $V$ be Hilbert spaces. Let $P$ be an operator such that
\[
P=\begin{bmatrix} P_{11} & P_{12} \\ P_{21} & P_{22} \end{bmatrix}:H\oplus U\to G\oplus V
\] 
and $||P|| <1$. Let $X:V\to U$ be an operator such that $||X||\leq1$ and $I-P_{22}X$ is invertible.
Then $||\mathcal{F}_{P}(X)|| <1$. 
\end{cor}

\begin{proof}
By Proposition \ref{calF},
\begin{align*}
I_H-\mathcal{F}_P(X)^*&\mathcal{F}_P(X)=\\
&=P^*_{21}(I_V-X^*P^*_{22})^{-1}(I_V-X^*X)(I_V-P_{22}X)^{-1}P_{21}\\
&\qquad+\begin{bmatrix} I_H&P^*_{21}(I_V-X^*P^*_{22})^{-1}X^* \end{bmatrix}(I_{H\oplus U}-P^*P)
\begin{bmatrix} I_H\\X(I_V-P_{22}X)^{-1}P_{21} \end{bmatrix}.
\end{align*}
Let $A=(I_V-P_{22}X)^{-1}P_{21}:H\to V$ and
\[
B=\begin{bmatrix} I_H\\X(I_V-P_{22}X)^{-1}P_{21} \end{bmatrix}=\begin{bmatrix} I_H\\XA \end{bmatrix}:H\to H\oplus U.
\]
Then 
\[
I_H-\mathcal{F}_P(X)^*\mathcal{F}_P(X)=A^*(I_V-X^*X)A+B^*(I_{H\oplus U}-P^*P)B.
\]

Since $||X||\leq1$, $B$ has zero kernel and $||P|| < 1$,  
we have 
\[
A^*(I_V-X^*X)A\geq0\text{ and }B^*(I_{H\oplus U}-P^*P)B > 0.
\] 
Thus  
\[
I_H-\mathcal{F}_P(X)^*\mathcal{F}_P(X)  >  0.
\] 
Therefore $||\mathcal{F}_P(X)|| < 1$.
\end{proof}

\begin{rem}\label{invertpartdefinedrealform}
\textup{Suppose, in addition, $||X||<1$ in Corollary \ref{calFleq1}. Then  
\[||P_{22}X||\leq||P_{22}||\, ||X||\leq||X||<1\] and so $I_V-P_{22}X$ is automatically invertible.} 
\end{rem}

\begin{rem}\label{calFholo} 
\textup{Suppose, in addition, $H=G=\Bbb{C}^n$, $U=V$ and $X=z\cdot I_V$ in Corollary \ref{calFleq1}. Then,  
by Remark \ref{invertpartdefinedrealform}, $I_V-zP_{22}$ is invertible for all $z\in\Bbb{D}$. Moreover, 
the linear fractional transformation $\mathcal{F}_P$, given by
\[\mathcal{F}_P(z)=P_{11}+zP_{12}(I_V-zP_{22})^{-1}P_{21}\text{ for all }z\in\Bbb{D},\] 
is holomorphic on $\Bbb{D}$.} 
\end{rem}

\end{document}